\documentclass[a4paper,12pt]{article} 
\usepackage{amsmath,amsfonts} %symboles mathematiques% 
\usepackage [french] {babel}
\usepackage {ucs}
\usepackage[utf8x]{inputenc}

\newtheorem{thm}{Theorem}[section]   
   
\newtheorem{prop}[thm]{Proposition}   
\newtheorem{lm}[thm]{Lemma}

\newcommand{\RR}{\mathbb{R}}   
\newcommand{\di}{\displaystyle}

\newcommand{\NN}{\mathbb{N}}

\renewcommand{\epsilon}{\varepsilon}
\newcommand{\e}{\varepsilon}

\begin{document}   
   
\title{\textbf{Nontrivial dynamics beyond the logarithmic shift}\\ 
\textbf{in two-dimensional Fisher-KPP equations}}

\author{
{\bf Jean-Michel Roquejoffre} \\ 
Institut de Math\'ematiques de Toulouse (UMR CNRS 5219) \\ 
Universit\'e Toulouse III,
118 route de Narbonne\\
31062 Toulouse cedex, France \\ 
\texttt{jean-michel.roquejoffre@math.univ-toulouse.fr}\\ 
\\[2mm]
{\bf Violaine Roussier-Michon} \\ 
Institut de Math\'ematiques de Toulouse (UMR CNRS 5219)\\ 
INSA Toulouse,
135 av. Rangueil \\
31077 Toulouse cedex 4, France\\ 
\texttt{roussier@insa-toulouse.fr}}
   
\maketitle   
\begin{abstract}
We study the asymptotic behaviour, as time goes to infinity, of the Fisher-KPP equation $\partial_t u=\Delta u +u-u^2$ in spatial dimension $2$, when the initial condition looks like a Heaviside function. Thus the solution is, asymptotically in time, trapped between two planar critical waves whose positions are corrected by the Bramson logarithmic shift. The issue is whether, in this reference frame, the solutions will converge to a travelling wave, or will exhibit more complex behaviours. We prove here that both convergence and nonconvergence may happen: the solution may converge towards one translate of the planar wave, or oscillate between two of its  translates. This relies on the behaviour of the initial condition at infinity in the transverse direction.
\end{abstract}

%%%%%%%%%%%%%%%%%%%%%%%%%%%%%%%%%%%%%%%%%%%%%%%%%%%%%%%%%%%%%%%%%%%%%%
%%%%%%%%%%%%%%%%%%%%%%%%%%%%%%%%%%%%%%%%%%%%%%%%%%%%%%%%%%%%%%%%%%%%%%
\section{Introduction}    
%%%%%%%%%%%%%%%%%%%%%%%%%%%%%%%%%%%%%%%%%%%%%%%%%%%%%%%%%%%%%%%%%%%%%
%%%%%%%%%%%%%%%%%%%%%%%%%%%%%%%%%%%%%%%%%%%%%%%%%%%%%%%%%%%%%%%%%%%%

\noindent
The paper is devoted to the large time behaviour of the solution of the reaction-diffusion equation 
\begin{align}
\label{KPP} \partial_t u= \Delta u +f(u), & \quad t>1 \, , \,  (x,y) \in \RR^2\\
\notag u(1,x,y)=u_0(x,y), & \quad \quad \quad  \quad (x,y) \in \RR^2
\end{align}
We will take
$$
f(u)=u(1-u) \,  \mbox{ if } u\in [0,1] \mbox{ and } f(u)=0 \mbox{ if } u \notin [0,1];
$$
thus  $f $ is said to be of the Fisher-KPP type.The initial datum $u_0$ is in ${\cal C}(\RR^2)$ and there exist $x_2<x_1$ such that
\begin{equation}
\label{assumption}
1-H(x-x_2)\leq u_0(x,y)\leq 1-H(x-x_1)
\end{equation}
where $H$ is the Heaviside function. 
Then, since $f$ is globally Lipschitz on $\mathbb{R}$, there exists a unique classical solution $u(t,x,y)$ in $ {\cal C}([1,+\infty[ \times \RR^2, (0,1))$ to equation \eqref{KPP} emanating from such $u_0$.

\noindent
The assumptions on $f$ imply that $0$ and $1$ are, respectively, unstable and stable equilibria for the ODE $\dot{\zeta}=f(\zeta)$. For the PDE \eqref{KPP}, the state $u\equiv 1$ invades the state $u \equiv 0$. 
Equation \eqref{KPP} admits one-dimensional travelling fronts $U(x-ct)$ if and only if $c\geq c^*=2$ where the profile $U$, depending on $c$, satisfies
\begin{equation}
\label{onde}
U'' +c \,U'+f(U)=0 , \quad x \in  \RR,
\end{equation}
together with the boundary conditions at infinity 
\begin{equation} 
\label{CL onde}
\lim\limits_{x \to - \infty} U(x)=1 \mbox{ and } 
\lim\limits_{x \to + \infty} U(x)=0. 
\end{equation} 
Any solution $U$ to \eqref{onde}-\eqref{CL onde} is a shift of a fixed profile $U_c$: $U(x)=U_c(x+s)$ with some fixed $s \in \RR$. The profile $U_{c^*}$ at minimal speed $c^*=2$ satisfies
$$
U_{c^*}(x)=(x+k) \,  e^{- x} +O(e^{-(1+\delta_0)x})\, , \, \mbox{ as } x \to +\infty
$$
for some universal constant $k \in \RR$, and some small $\delta_0>0$.

%%%%%%%%%%%%%%%%%%%%%%%%%%%%%%%%%%%%%%%%%%%%%%%%%%
\subsection{Convergence for the 1D KPP equation: related works}
%%%%%%%%%%%%%%%%%%%%%%%%%%%%%%%%%%%%%%%%%%%%%%%%%%

\noindent
The large time behaviour of the one dimensional  problem 
\begin{equation}
\label{KPP1D}
\partial_t u=\partial_{xx} u +f(u),\quad t>1\, , \, x\in\RR
\end{equation}
has a history of important contributions. One of the 
first, and perhaps most well-known one, is 
 the pioneering KPP paper \cite{KPP}. Kolmogorov, Petrovskii and Piskunov proved 
that the solution of~\eqref{KPP1D}, starting from $1-H(x)$, converges to $U_{c*}$ in shape: there is a   function
$$
\sigma_\infty(t)=2t+o(t),
$$
such that 
\[
\di
\lim_{t\to+\infty}u(t,x + \sigma_\infty(t))=U_{c^*}(x) \quad\mbox{  uniformly in }x \in \RR.
\]
The main ingredient in \cite{KPP} is the monotonicity of $\partial_x u$ on the level sets of $u$. This 
argument was recently revisited by Ducrot-Giletti-Matano \cite{DGM}, Nadin \cite{N}, 
for results in the same spirit, concerning one-dimensional inhomogeneous models.

The second one makes precise the $\sigma_\infty(t)$: in \cite{Bramson1,Bramson2},   Bramson
 proves the following
\begin{thm}
\label{thm bramson}
There is  a constant $x_\infty$, depending on $u_0$, such that
\[
\sigma_\infty(t)=2t-\di
\frac32\ln t-x_\infty+o(1),\hbox{ as } t\to+\infty.
\]
\end{thm}
Theorem \ref{thm bramson} was proved through elaborate probabilistic arguments.  
A natural question was thus to prove Theorem \ref{thm bramson} with purely PDE arguments. 
In that spirit, a weaker version,
precise up to the~$O(1)$ term,  is the main result of~\cite{HNRR1} (which is actually the PDE counterpart of \cite{Bramson1}):
$$\sigma(t)=2t-\di\frac32\ln t+O(1) \hbox{ as } t\to+\infty \, .$$
This was extended   for the much more
difficult case of the periodic in space
coefficients, see \cite{HNRR2}. Bramson's theorem \ref{thm bramson} is fully recovered in 
\cite{NRR-Brezis}, with once again  simple and robust PDE arguments. The dynamics beyond the shift has also been the subject of 
intense studies: let us mention the paper \cite{ES}, which proposes a universal behaviour for $\sigma(t)-\sigma_\infty(t)$, by means of formal asymptotic arguments.
See also \cite{VS}. The universal correction is obtained, in a mathematically rigorous way, in \cite{NRR2}. See also \cite{BBHR} for asymptotics in a related free boundary problem.

%%%%%%%%%%%%%%%%%%%%%%%%%%%%%%%%%%%%%%%%%%%%%%%%%%%
\subsection{Question and results}
%%%%%%%%%%%%%%%%%%%%%%%%%%%%%%%%%%%%%%%%%%%%%%%%%%

\noindent
Let us come back to our two-dimensional case.  Let $u_i(t,x)$, $i\in\{1,2\}$ be the solution of the one-dimensional problem \eqref{KPP1D} emanating from $u_i(1,x)=1-H(x-x_i)$. By the maximum principle we have $u_2(t,x)\leq u(t,x,y)\leq u_1(t,x)$. And so, there exist $x_{\infty,1}\geq x_{\infty,2}$ such that, if  an arbitrary level set of $u(t,.)$ is represented by the graph $\{x=\sigma(t,y)\}$ - this is not always true, but certainly true if $u_0$ is nonincreasing in $x$ (applying the maximum principle on $u_x$)- there is a   function $\sigma_\infty(t,y)\in [x_{2,\infty},x_{1,\infty}]$ such that
\begin{equation}
\label{sigma}
\sigma(t,y)=2t-\frac32 \ln t+\sigma_\infty(t,y).
\end{equation}
The issue is: does this function $\sigma_\infty$ converge for large times? In one space dimension ($\sigma_\infty$ only depending on time), this is true. In order to realise that it is a true issue in two space dimensions, let us make a parallel with the case where $f$ is bistable: there is $\theta\in(0,1)$ such that $f(u)<0$ if $u\in(0,\theta)$ and $f(u)>0$ on $(\theta,1)$. Contrary to the KPP case, the travelling wave problem \eqref{onde}-\eqref{CL onde} has a unique orbit $(c_*,U_{c_*})$. If $u(1,x)=1-H(x)$, then (Fife-McLeod \cite{FML}) $u(t,x)$ converges exponentially fast to the wave profile; in other words there are $x_\infty\in\RR$ and $\omega>0$ such that
$$
u(t,x)=U_{c_*}(x-c_*t +x_\infty)+O(e^{-\omega t}) \mbox{ uniformly in } x \in \RR.
$$
However, under the assumption \eqref{assumption}, and if $\sigma(t,y)$ denotes any level set of $u(t,.)$, there is (Roquejoffre, Roussier-Michon \cite{JMR-VRM}) a bounded function $\sigma_\infty(t,y)$ such that 
$$\sigma(t,y)=c_*t+\sigma_\infty(t,y)+O(t^{-1/2}),$$
and, depending on the initial datum $u_0$, the function $\sigma_\infty(t,y)$ may or not converge as time goes to infinity. It is therefore legitimate to suspect a phenomenon of that kind here, and this is exactly what happens.

Let us now state and explain our results. The first one says that the large time dynamics of \eqref{KPP} is, in some sense, that of the heat equation.
%%%%%%%%%%%%%%%%%%%%%%%%%%%%%%%%%%%%%
\begin{thm}
\label{thm1}
Let $u_0$ satisfy assumption \eqref{assumption}. For every small $\e>0$, there is $T_\e>0$ and a function $a_0^\e$,  with $\Vert a_0^\e\Vert_\infty$ and $\Vert\di da_0^\e/dy\Vert_\infty$ bounded in $\e$, such that the solution $u$ of \eqref{KPP} emanating from $u_0$ satisfies
$$u(t,x,y)=U_{c_*}\biggl(x-2t+\frac32{\mathrm{ln}}t-{\mathrm{ln}}(a^\e(t,y)+O(\e))\biggl)+O(t^{-1/2}),\quad\hbox{for $t\geq T_\e$,}$$
where the function $a^\e(t,y)$ solves the heat equation
$$(\partial_t-\partial_{yy}) a^\e=0\, , \quad t>1 \, ,\, y\in\RR \, , \quad\quad a^\e(1,y)=a_0^\e(y).$$
\end{thm}
%%%%%%%%%%%%%%%%%%%%%%%%%%%%%%%%%%%%%%

\noindent
This explains that \eqref{KPP} has, beyond the logarithmic shift, a large time dynamics which mimics that of the heat equation. We point out that this result is
optimal, since the solution of the heat equation does not, in general, converge to anything: see for instance Collet-Eckmann \cite{CE}, V\`azquez-Zuazua \cite{VZ}.
We will, by the way, use those results to construct solutions that do not converge beyond the shift.

Theorem \ref{thm1} is the most general one can prove. However, it does not really say whether the solution will, or not, converge to something, for the simple reason that it does not exclude a sequence $(a_0^\e)_\e$ such that the heat equation starting from $a_0^\e$ will diverge for $\e=O(1)$, and converge to something as $\e$ becomes very small. So, in the following result, we are going to show that both types of behaviour may happen: convergence to a single wave, or, on the contrary,   nonconvergence. Let us not forget, though, that the asymptotic dynamics is that of the heat equation. So, nonconvergence will occur through infinitely slow oscillations between two waves. Assume, for definiteness, that $u_0$ is nonincreasing in $x$. This is by no means necessary but, since we are not aiming for utmost generality, this slight loss of generality will be compensated by a lighter formulation. Let $\sigma_\infty(t,y)$ be given by \eqref{sigma}.
\begin{thm}
\label{thm2}
The following situations hold.
\begin{enumerate}
\item There are initial data $u_0(x,y)$, satisfying  assumptions \eqref{assumption}, such that $t\mapsto\sigma_\infty(t,0)$ does not converge as $t\to+\infty$. 
\item Assume the existence of two functions $u_0^\pm(x)$, and $x_1\leq x_2$, such that 
$$1-H(x-x_1)\leq u_0^+(x),\ u_0^-(x)\leq 1-H(x-x_2),$$
and such that
$$
\lim_{y\to\pm\infty}u_0(x,y)=u_0^\pm(x),\ \hbox{uniformly in $x\in\RR$.}
$$
Let the constants $\sigma^\pm_\infty$ be defined as follows: if $u^\pm(t,x)$ is the solution of \eqref{KPP1D} emanating from $u_0^\pm(t,x)$, then
$$u^\pm(t,x)=U_{c_*}\biggl(x-2t+\frac32{\mathrm{ln}}t+\sigma_\infty^\pm\biggl)+o_{t\to+\infty}(1).$$
Then we have
$$
\lim_{t\to+\infty}\sigma_\infty(t,y)=-{\mathrm{ln}}\biggl(\frac{e^{-\sigma_\infty^+}+e^{-\sigma_\infty^-}}2\biggl),
$$
uniformly on every compact set in $y$. If $\sigma^+_\infty=\sigma^-_\infty$, the convergence is uniform in $y$.
\item Assume the existence of $u_\infty(x,y)$, periodic in $y$, such that 
$$
\lim_{y\to +\infty}\biggl(u_0(x,y)-u_\infty(x,y)\biggl)=0,\ \hbox{uniformly in $x$.}
$$
Then $\sigma_\infty(t,y)$ converges to a constant as $t\to+\infty$, uniformly in $y$.
\end{enumerate}
\end{thm}
We could of course imagine more situations, such as, for instance, the existence of two periodic functions $u^\pm_\infty(x,y)$ such that  $u_0(x,y)$ resembles $u^+_\infty(x,y)$ (resp. $u^-_\infty(x,y)$) as $y\to+\infty$ (resp. $y\to-\infty$)...
Another interesting question is to understand what happens beyond $\sigma_\infty(t,y)$, in other words can one devise an asymptotic expansion, which could hold only uniformly on every compact in $y$.

%%%%%%%%%%%%%%%%%%%%%%%%%%%%%%%%%%%%%%%%%%%%%%%%%%%
\subsection{Other multi-D configurations}
%%%%%%%%%%%%%%%%%%%%%%%%%%%%%%%%%%%%%%%%%%%%%%%%%%%

\noindent
Let us briefly mention the state of the art when the initial data, instead of being trapped between two transates of the Heaviside function, is compactly suppported,
and let us restrict ourselves to \eqref{KPP} - we do not assume the medium to be heterogeneous. The first, and most general result, is due to Aronson-Weinberger \cite{aw}. The solution $u$ spreads at the speed $c^*=2\sqrt{f'(0)}=2$ in the sense that 
$$
\min_{|x|\leq ct} u(t,x) \to 1 \mbox{ as } t \to +\infty \, , \mbox{ for all } 0 \leq c < c^*
$$
and 
$$
\sup_{x\geq ct} u(t,x) \to 0 \mbox{ as } t \to +\infty \, , \mbox{ for all }  c > c^*
$$
This estimate is made precise  up to $O(1)$ terms in G\"artner \cite{Ga}. See also Ducrot \cite{D}, who uses the ideas of \cite{HNRR1} to give a PDE proof of \cite{Ga}. In fact, the precise large-time behaviour in the bistable case is known (Roussier-Michon \cite{Roussier}, Yagisita \cite{yagisita}). The extension of these results to the KPP case in underway \cite{RRR}.

\noindent
{\bf Acknowledgements.} JMR  is supported by  the European Union's Seventh 
Framework Programme (FP/2007-2013) / ERC Grant
Agreement n. 321186 - ReaDi - ``Reaction-Diffusion Equations, Propagation and Modelling''. Both authors are supported by the ANR project NONLOCAL ANR-14-CE25-0013.

%%%%%%%%%%%%%%%%%%%%%%%%%%%%%%%%%%%%%%%%%%%%%%%%%%%
%%%%%%%%%%%%%%%%%%%%%%%%%%%%%%%%%%%%%%%%%%%%%%%%%%%
\section{Equations, strategy of the proof, organisation of the paper}
\label{section equations}
%%%%%%%%%%%%%%%%%%%%%%%%%%%%%%%%%%%%%%%%%%%%%%%%%%%
%%%%%%%%%%%%%%%%%%%%%%%%%%%%%%%%%%%%%%%%%%%%%%%%%%%

\noindent
There is a sequence of transformations that bring the equations under the \eqref{KPP} to 
a form that will be amenable to treatment. 
\begin{enumerate}
\item We observe the equation \eqref{KPP} in the reference frame whose origin is $X(t)=2t- \frac{3}{2} \ln t $ and choose the change of variables $x'=x-X(t)$ and $u(t,x,y)=u_1(t,x-X(t),y)$. After dropping the primes and indexes, equation \eqref{KPP} becomes
\begin{equation}
\label{ref mvt} 
\partial_t u= \Delta u + \left( 2-\frac{3}{2t}\right) \partial_x u + u-u^2 \, , \quad t>1  \, ,  \quad (x,y) \in \RR^2
\end{equation}
with initial datum $u(1,x,y)=u_0(x+2,y)$.
 \item To follow the exponential decrease of the wave $U_{c^*}$, it will be useful to take it out and set $u(t,x,y)=e^{-x} v(t,x,y)$; \eqref{ref mvt} thus becomes
\begin{equation}
\label{dec onde} 
\partial_t v= \Delta v - \frac{3}{2t} \left(\partial_x v - v\right) -e^{-x} v^2  \, , \quad t>1  \, ,  \quad (x,y) \in \RR^2
\end{equation}
with initial datum $v(1,x,y)=e^x u_0(x+2,y)$.
\item Finally, if we want to study \eqref{dec onde} in the diffusive zone, i.e. the region $x\sim\sqrt t$,
 we introduce self similar variables 
$\displaystyle \xi=\frac{x}{\sqrt{t}}$, $\tau=\ln t $. The variable $y$ is unchanged:
\begin{equation}
\label{var auto sim}
w(\tau,\xi,y)=w \left(\ln t, \frac{x}{\sqrt{t}}, y \right)=\frac{1}{\sqrt{t}} v(t,x,y)
\end{equation}
Then \eqref{dec onde} becomes 
\begin{equation}
\label{shift} 
 \partial_{\tau} w = {\cal L} w + e^{\tau} \partial_{yy} w -\frac{3}{2} e^{-\frac{\tau}{2}} \partial_{\xi} w - e^{\frac{3}{2}\tau - \xi e^{\frac{\tau}{2}}} w^2  \, , \quad \tau>0 \, ,  \quad (\xi,y) \in \RR^2
\end{equation}
where 
$${\cal L} w = \partial_{\xi\xi} w +\frac{\xi}{2} \partial_{\xi} w + w$$
with initial datum $w(0,\xi,y)=e^{\xi} u_0(\xi+2,y)$.
\end{enumerate}
In the sequel, we will use the form that will be best suited to our purposes.
Let us say a word about the strategy of the proof of Theorem \ref{thm1}. In one space dimension, \eqref{shift} becomes
 $$\partial_{\tau} w = {\cal L} w -\frac{3}{2} e^{-\frac{\tau}{2}} \partial_{\xi} w - e^{\frac{3}{2}\tau - \xi e^{\frac{\tau}{2}}} w^2 \, , \quad  \tau>0 \, ,\quad \xi\in \RR.$$
The main step of the proof in \cite{NRR-Brezis} was to prove the existence of a constant $\alpha_\infty>0$ such that
$$
w(\tau,\xi)\longrightarrow_{\tau\to+\infty} \alpha_\infty\xi^+e^{-\xi^2/4},\quad\hbox{in $\{\xi\geq e^{-(\frac12-\delta)\tau}\}$},
$$
where $\delta>0$ is arbitrarily small. We would then define the translation $\sigma_\infty(t)$ such that 
$$
U_{c_*}(x+\sigma_\infty(t))\biggl\vert_{x=t^\delta}=e^{-x}v(t,x)\biggl\vert_{x=t^\delta}.
$$
That is, 
\begin{equation}
\label{alpha infini}
\sigma_\infty(t)=-{\mathrm{ln}}\alpha_\infty+O(t^{-\delta}).
\end{equation}
We would then prove the uniform convergence to $U_{c_*}(x-{\mathrm{ln}}\alpha_\infty)$ by examining the difference
$$
\tilde v(t,x)=\big\vert v(t,x)-U_{c_*}(x+\sigma_\infty(t))\big\vert
$$
in the region $\{x<t^\delta\}$. It turned out that $\tilde v(t,x)$ was a subsolution of (a perturbation of) the
heat equation
\begin{equation}
\label{e2.301}
\begin{array}{rll}
V_t= &V_{xx} + O(t^{1-\delta}) \, ,& \quad t>0 \, ,\, -t^\delta<x<t^\delta\\
V(t,-t^\delta)=&e^{-t^\delta} \, ,& \quad t>0\\
V(t,t^\delta)=&0 \, ,& \quad t>0.
\end{array}
\end{equation}
The condition on the left simply comes from the fact that $v(t,x)$ decays, by definition, like $e^x$ at $-\infty$. Although the domain looks very large, its first Dirichlet eigenvalue is of the order $t^{-2\delta}$, hence a much larger quantity than the right hand side of \eqref{e2.301}. Thus $V(t,x)$ goes to 0 uniformly in $x$ as $t\to+\infty$, which implies the sought for convergence result.

In what follows, we are going to adapt these ideas to our setting. The main additional difficulty is the transverse diffusion, which, in a very paradoxical way, does not help us. This is not a rhetorical argument: its presence is really what prevents convergence, in most cases. This implies that we will have to be quite careful with the estimates.

\medskip
The paper is organised as follows. In Section 3, we explain how the behaviour of $u(t,x,y)$ in the half plane $\{x<t^\delta,y\in\RR\}$ is slaved to that on the line $\{x=t^\delta,y\in\RR\}$. In Section 4, we characterise the asymptotic behaviour of a general linear equation that encompasses, in particular, equation \eqref{shift}. In Section 5, we define sub and super solutions that will enable us to prove Theorem \ref{thm1}. Finally, Theorem \ref{thm2} is proved in Section 6.

%%%%%%%%%%%%%%%%%%%%%%%%%%%%%%%%%%%%%%%%%%%%%%%%%%
%%%%%%%%%%%%%%%%%%%%%%%%%%%%%%%%%%%%%%%%%%%%%%%%%%
\section{Control of the solution by its value at $t^\delta$}
\label{droite}
%%%%%%%%%%%%%%%%%%%%%%%%%%%%%%%%%%%%%%%%%%%%%%%%%%%
%%%%%%%%%%%%%%%%%%%%%%%%%%%%%%%%%%%%%%%%%%%%%%%%%%%

\noindent
 The goal of this section is to prove, as announced in the introduction, that controlling the solution slightly to the right of the $O(1)$ in $x$ area implies, provided that the  control is well-tailored, the control of the solution to the entire region to the left. From now on we consider $\delta\in (0,\di\frac12)$, that will be as small as we wish.

 %%%%%%%%%%%%%%%%%%%%%%%%%%%%%%%%%%%%%%%%%%%%%%%%%%%%%
\subsection{The basic result}
%%%%%%%%%%%%%%%%%%%%%%%%%%%%%%%%%%%%%%%%%%%%%%%%%%%%

\noindent
Let $a(t,y)$  be a smooth function such that
\begin{itemize}
\item There are constants $0<\underline a_0 \leq\overline a_0<+\infty$ that bound $a$:
 \begin{equation}
 \label{hyp 1 sur a}
 \forall t>1 \, , \, \forall y \in \RR \, , \quad \underline a_0\leq a(t,y)\leq\overline a_0,
 \end{equation}
\item There is a constant $C_0>0$ depending on $\underline a_0$ and $\overline a_0$  that bounds the derivatives of $a$:
\begin{equation}
\label{hyp 2 sur a}
\forall t>1 \, , \, \forall y \in \RR \, , \quad
\vert\partial_y a(t,y)\vert\leq\frac{C_0}{\sqrt t},\quad \max(\vert\partial_{yy} a(t,y)\vert, \vert\partial_{t} a(t,y)\vert)\leq\frac{C_0}{t}.
\end{equation}
\end{itemize}
We define $\gamma(t,y)$ by the relation
\begin{equation}
\label{e3.1}
U_{c_*}(t^\delta+\gamma(t,y))=t^\delta e^{-t^\delta-1/4t^{1-2\delta}}\frac{a(t,y)}{\sqrt{2\sqrt{\pi}}}:=u_+^a(t,y).
\end{equation}
We have therefore, for large $t$ and $\delta \in (0,\frac{1}{3})$:
$$
\gamma(t,y)=-\ln \left( \frac{a(t,y)}{\sqrt{2\sqrt{\pi}}} \right)+O(t^{-\delta}).
$$
More important we have, from the implicit functions theorem, that 
$\gamma$ is at least $C^1$ in $t$ and $C^2$ in $y$, and we have, for a
universal constant $C$:
\begin{equation}
\label{e3.2}
\begin{array}{rll}
\vert\partial_y \gamma(t,y)\vert\leq &C\vert\partial_y a(t,y)\vert\\
\vert\partial_{yy}\gamma(t,y)\vert\leq&C\biggl(\vert\partial_{yy} a(t,y)\vert+(\partial_y a(t,y))^2\biggl)\\
\vert\partial_t\gamma(t,y)\vert\leq&C\biggl(\vert\partial_t a(t,y)\vert+\displaystyle\frac{a(t,y)}{t^{1-\delta}}\biggl)
\end{array}
\end{equation}
Let $u_a(t,x,y)$ be a solution of
\begin{align}
\label{e3.11}
\partial_t u_a=\Delta u_a +\left(2-\di\frac3{2t}\right) \partial_x u_a + u_a-u_a^2 & \quad t>1 \, , \, x\leq t^\delta \, , \, y\in\RR\\
\notag  u_a(t,t^\delta,y)=u_+^a(t,y) & \quad t\geq 1 \, ,  \, x=  t^\delta \, , \, y\in\RR\\
\notag \di\inf_{y\in\RR}\liminf_{x\to-\infty}u_a(1,x,y)>0 &
\end{align}
Here is the main result of this section. 
%%%%%%%%%%%%%%%%%%%%%%%%%%%%
\begin{thm}
\label{t3.1}
For $\delta \in (0,\frac14)$ and $u_a$ solution to equation \eqref{e3.11} where $u_+^a$ is defined in \eqref{e3.1} and $a$ satisfies assumptions \eqref{hyp 1 sur a} and \eqref{hyp 2 sur a}, we have for any $t>1$
$$
\di\sup_{|x|\leq t^\delta} \sup_{y \in \RR} e^x\biggl\vert u_a(t,x,y)-U_{c_*}(x+\gamma(t,y))\biggl\vert\leq \frac{C}{t^\lambda},
$$
for some universal constant $C>0$ and $\lambda \in (0, 1-4\delta)$.
\end{thm}
%%%%%%%%%%%%%%%%%%%%%%%%%%%
\noindent\textbf{ Proof.} We simply set
$$
s(t,x,y)=e^x\big( u_a(t,x,y)-U_{c_*}(x+\gamma(t,y))\big),
$$
Then,  for any $ t>1 $, $x<t^\delta$, and  $y \in \RR$
$$
\partial_t s-\Delta s +\frac{3}{2t} (\partial_x s-s) +s(u_a+U_{c_*}(x+\gamma))= e^x \left( (\partial_{yy} \gamma - \partial_t \gamma)U' + (\partial_y \gamma)^2 U'' \right)
$$
so that by \eqref{e3.2}, we have
\begin{align}
\label{e3.4}
\partial_t s-\Delta s +\frac{3}{2t} (\partial_x s-s) +s(u_a+U_{c_*}(x+\gamma))= O\left(\frac{1}{t^{1-2\delta}} \right) & \quad t>1 \, , \, x<t^\delta \, , \, y \in \RR \\
\notag s(t,t^\delta,y)=0 & \quad t>1 \, , \, x=t^\delta \, , \, y \in \RR \\
\notag \di\sup_{y \in \RR} s(t,-t^{\delta},y)=O\left(e^{-t^\delta} \right) & \quad t>1 \, , \, x=-t^\delta \, , \, y \in \RR
\end{align}
The last equation comes from the definition of $s$, as the product of a bounded function by an exponential.
As in \cite{NRR-Brezis}, a super-solution to \eqref{e3.4} is devised as
$$
\overline s(t,x,y)=\frac{A}{t^\lambda}  \cos \left(\frac{x}{t^{\delta+\tilde\e}}\right),
$$
where $\delta \in (0,\frac14)$, $\lambda \in(0,1-4\delta)$,  $\tilde\e >0 $ is small enough such that $2\delta +2\tilde{\e} +1- \lambda <1-2\delta$ and $A>0$ large enough. The idea is that the first Dirichlet eigenvalue of $(-\partial_{xx})$ in the interval $(-t^\delta,t^\delta)$ is of order $t^{-2\delta}$ (a nonintegrable power of $t$ if $\delta$ is small enough), whereas the right hand side of \eqref{e3.4} is of the order $t^{2\delta-1}$, a much larger power. And so, $\overline s$ will dominate $s$, which proves the result. \rule{2mm}{2mm}

%%%%%%%%%%%%%%%%%%%%%%%%%%%%%%%%%%%%%%%%%%%%%%%%
\subsection{Perturbative results}
%%%%%%%%%%%%%%%%%%%%%%%%%%%%%%%%%%%%%%%%%%%%%%%%

\noindent
Consider $\e>0$ and $b(t,y)$ a smooth function such that for any $t>1$ and $y \in \RR$
\begin{equation}
\label{e3.10}
\vert b(t,y)\vert\leq\e+\frac{C}{t^\delta},
\end{equation}
for some constant $C>0$. Note that no assumption is made on the derivatives of $b$ and, in particular, no assumption on a possible time decay of $\partial_t b$ or $\partial_y b$. Set, this time
\begin{equation}
\label{e3.30}
u_+^{a+b}(t,y)=t^\delta e^{-t^\delta-1/4t^{1-2\delta}}\frac{a(t,y)+b(t,y)}{\sqrt{2\sqrt\pi}}.
\end{equation}
Theorem \ref{t3.1} perturbs into the following 
\begin{prop}
\label{p3.1}
For $\delta \in (0, \frac15)$, let $u_a$ (resp. $u_{a+b}$) be a solution of the Dirichlet problem \eqref{e3.11}, with  boundary condition $u_+^a(t,y)$ (resp. $u_+^{a+b}$).  There exists $C>0$, depending on $u_a(1,.)$ and $u_{a+b}(1,.)$ such that
for any $t>1$ 
$$
\di\sup_{|x|\leq t^\delta} \sup_{y \in \RR} e^x \vert u_{a+b}(t,x,y)- u_a(t,x,y)\vert\leq C(\e+\frac1{t^\delta}).
$$
\end{prop}
{\bf Proof.} Define $\underline u(t,x,y)$ (resp. $\overline u(t,x,y)$) as the solutions of \eqref{e3.11} with the following data:
$$\left\{
\begin{array}{rll}
\underline u(t,t^\delta,y)=& u_+^{a+b} -C(\e+t^{-\delta}),\quad \underline u(1,x,y)=\min\biggl(u_a(1,x,y), u_{a+b}(1,x,y)\biggl)\\
\overline u(t,t^\delta,y)=& u_+^{a+b} +C(\e+t^{-\delta}),\quad \overline u(1,x,y)=\max\biggl(u_a(1,x,y), u_{a+b}(1,x,y)\biggl)\\
\end{array}
\right.$$
Both $\overline u$ and $\underline u$ fall in the assumptions of Theorem \ref{t3.1}, thus $\overline u$ approaches $U_{c_*}(x+\overline\gamma(t,y))$
(resp. $\underline u$ approaches $U_{c_*}(x+\underline\gamma(t,y))$ like $t^{-\lambda}$ as $t\to+\infty$ with $\lambda \in (0,1-4\delta)$. The definition of $\overline \gamma$ and $\underline 
\gamma$ mimick that of $\gamma$ in the preceding section; in other words the translation of $U_{c_*}$ is adjusted to coincide with the solution at the boundary.
Thus we have
$$
\vert \overline\gamma(t,y)-\underline\gamma(t,y)\vert\leq C(\e+t^{-\delta}),
$$
and the proposition follows since $1-4\delta >\delta$. \rule{2mm}{2mm}

%%%%%%%%%%%%%%%%%%%%%%%%%%%%%%%%%%%%%%%%%%%%%%%%%%%%
%%%%%%%%%%%%%%%%%%%%%%%%%%%%%%%%%%%%%%%%%%%%%%%%%%%%%
\section{A Dirichlet problem in the diffusive zone}
%%%%%%%%%%%%%%%%%%%%%%%%%%%%%%%%%%%%%%%%%%%%%%%%%%%%
%%%%%%%%%%%%%%%%%%%%%%%%%%%%%%%%%%%%%%%%%%%%%%%%%%%

\noindent
Consider the following equation for $\e >0$, $\lambda>0$, 
\begin{align}
\label{dirichlet} 
 \partial_{\tau} v = {\cal L} v + \frac{e^{\tau}}{\e^2} \partial_{yy} v + \e^{2\lambda} e^{-\lambda\tau} \left( \phi_{\e}(\tau) v + \psi_\e(\tau) \partial_{\xi} v +f_\e (\tau,\xi) \right) &\, ,  \quad \tau>0 \, , \, \xi >0 \, , \,  y \in \RR \\
 \notag v(\tau,0,y)=0 &\, , \quad  \tau >0 \, , \,  \xi=0 \, , \, y \in \RR  \\
 \notag v(0,\xi,y)=v_0(\xi,y) &\, , \quad\tau =0 \, , \,  \xi >0 \, , \,  y \in \RR 
\end{align}

%%%%%%%%%%%%%%%%%%%%%%%%%%%%%%%%%%%%%%%%%%%%%%%%%
\subsection{Behaviour for general initial data}
%%%%%%%%%%%%%%%%%%%%%%%%%%%%%%%%%%%%%%%%%%%%%%%%%

\noindent
With no particular assumption on the behaviour of $v_0$ in the direction $y$, we are going to prove the following result.
%%%%%%%%%%%%%%%%%%%%%%%%%%%%%%%%
\begin{thm}
\label{lm 2.2}
For any $\lambda>0$ and any $C_0>0$, there exist $\epsilon_0>0$ and $C>0$ such that  for any compact set $K \subset \RR^+$, there is $C_K>0$ such that  for any $\epsilon \in (0,\epsilon_0)$, any initial data $v_0 \in X$, any functions $\phi_\e$ and  $\psi_\e$ uniformly bounded in $\tau$ and $\epsilon$ by $C_0$ and any function $f_\e$ compactly supported in $\xi$ and uniformly bounded in $\e$ by $C_0$, there exists a unique solution $v \in {\cal C}(\RR^+, X)$ to \eqref{dirichlet} emanating from $v_0$ which satisfies 
$$
\forall \tau > 0 \, , \, \xi >0 \, , \, y \in \RR \, , \quad v(\tau,\xi,y)= \xi \left( \frac{e^{-\xi^2/4}}{\sqrt{2\sqrt{\pi}}} \left(\alpha_c(\tau,y) +\beta(\tau,y)\right)+ e^{- \frac{\lambda}{2} \tau} \tilde{v}(\tau,\xi,y) \right)
 $$
where  for any $\tau >0$,  $y\in\RR$
\begin{align*}
&\partial_{\tau} \alpha_c = \frac{e^{\tau}}{\e^2} \partial_{yy} \alpha_c \, , \quad \alpha_c(0,y)= \frac{1}{\sqrt{2\sqrt{\pi}}}  \int_0^{+\infty}\xi  v_0(\xi,y) d\xi\\
& \|\beta(\tau)\|_{L^\infty(\RR)} \leq C \epsilon^{2\lambda} \, \, \quad \qquad \qquad \|\partial_\tau \beta(\tau)\|_{L^\infty} \leq C \epsilon^{2\lambda} \\
& \|\partial_y \beta(\tau)\|_{L^\infty(\RR)} \leq C \epsilon^{2\lambda+1}  e^{-\frac{\tau}{2}} \qquad  \|\partial_{yy}\beta (\tau)\|_{L^\infty} \leq C \epsilon^{2\lambda+2} e^{-\tau}   
\end{align*}
and for any $\tau >0$, $\xi \in K$, $y\in\RR$
\begin{align*}
\mbox{max }(|\tilde{v}(\tau,\xi,y)|, |\partial_\tau \tilde{v}| \, ,& \, |\partial_\xi \tilde{v} | \,  , \, |\partial_{\xi\xi} \tilde{v} | ) \leq C_K  \e^{\lambda}\\
|\partial_y \tilde{v}| \leq C_K \e^{\lambda+1} e^{-\frac{\tau}{2}} \, , &\quad   |\partial_{yy} \tilde{v}| \leq C_K \e^{\lambda+2} e^{-\tau} 
\end{align*}
\end{thm}
%%%%%%%%%%%%%%%%%%%%%%%%%%%%%%%%%%%

\noindent
\textbf{Proof of theorem \ref{lm 2.2}. }Choose $\lambda >0$ and $C_0>0$. Set $\epsilon >0$ and consider $\phi_\e$, $\psi_\e$ and $f_\e$ uniformly bounded in $\tau$ and $\e$ by $C_0$. Assume also $f_\e$ is compactly supported in $\xi$.
Let $v$ be the solution to \eqref{dirichlet} emanating from $v_0 \in X$. Let us introduce the new function $w(\tau,\xi,y)=e^{\frac{\xi^2}{8}}v(\tau,\xi,y)$. This new function solves for any $\tau >0$, $\xi >0$ and $y \in \RR$.
\begin{align}
\label{dirichlet sym} 
 \partial_{\tau} w = {\cal M} w + \frac{e^{\tau}}{\e^2} \partial_{yy} w + \e^{2\lambda} \, e^{-\lambda \tau} & \left((\phi_\e(\tau)-\frac{\xi}{4}\psi_\e(\tau))w+ \psi_\e(\tau) \partial_\xi w + e^{\frac{\xi^2}{8}}f_\e(\tau,\xi)\right)   \\ 
\notag
w(\tau,0,y)=0 & \quad   \tau>0 \, , \, \xi=0 \, , \, y \in \RR  \\
\notag w(0,\xi,y)=w_0(\xi, y)=e^{\frac{\xi^2}{8}} v_0(\xi,y) & \quad \tau =0\, , \,  \xi >0 \, , \, y \in \RR 
\end{align}
where  ${\cal M} w= \partial_{\xi\xi} w +\left(\frac{3}{4}-\frac{\xi^2}{16}\right)  w$. Thus ${\cal D}({\cal M})=\{w \in H_0^2(\RR^+) \, | \, \xi^2 w \in L^2(\RR^+) \}$, $M$ is symmetric and its null space is generated by the unit eigenfunction $e_0(\xi)=\frac{1}{\sqrt{2\sqrt{\pi}}} \xi e^{-\frac{\xi^2}{8}}$. 
This linear operator defines a quadratic form on $ \{w \in H_0^1(\RR^+) \, | \, \xi^2 w \in L^2(\RR^+) \}$ as
$$
q(w)=<-{\cal M}w,w>_{L^2(\RR^+)}=\int_0^{+\infty} (\partial_{\xi} w)^2 + \left(\frac{\xi^2}{16} - \frac{3}{4} \right) w^2 \, d\xi
$$
which is nonnegative and satisfies
$$
q(w) \geq \|w\|^2_{L^2(\RR^+)} \mbox{ if } <w,e_0>_{L^2(\RR^+)}=0
$$

%%%%%%%%%%%%%%%%%%%%%%%%%%
\begin{lm}
\label{w bounded}
There exist $\epsilon_0>0$ (depending on $\lambda$ and  $C_0)$ and  $C>0$ such that for any $\epsilon \in (0, \epsilon_0)$,  any $w$ solution to \eqref{dirichlet sym} emanating from $w_0 \in L^2(\RR^+,L^{\infty}(\RR))$ satisfies
$$
\forall \tau \geq 0\, , \quad \|w(\tau)\|_{L^2(\RR^+,L^{\infty}(\RR))} \leq C (\|w_0\|_{L^2(\RR^+,L^\infty(\RR))} + \e^{\lambda})
$$
\end{lm}
%%%%%%%%%%%%%%%%%%%%%%%%%%

\noindent
\textbf{Proof of lemma \ref{w bounded}.} Taking the $L^2(\RR^+)$ scalar product of \eqref{dirichlet sym} with $w$ leads to
\begin{align*}
\partial_{\tau} \|w\|_{L^2(\RR^+)}^2 + 2 q(w)= &  \frac{e^{\tau}}{\e^2}\left( \partial_{yy}\|w\|_{L^2(\RR^+)}^2 - 2 \|\partial_y w\|_{L^2(\RR^+)}^2\right) \\
& +  2 \epsilon^{2\lambda} e^{-\lambda\tau} \left( \phi_\e(\tau) \|w\|_{L^2(\RR^+)}^2 -\psi_\e(\tau) \int_0^{+\infty} \frac{\xi}{4} w^2 d\xi + \int_0^\infty e^{\frac{\xi^2}{8}} f_\e(\tau,\xi) w \, d\xi \right)
\end{align*}
Note that
$$
\int_0^{+\infty} \frac{\xi}{4} w^2 d\xi \leq \int_0^{+\infty} \left(\frac{\xi^2}{16}+\frac{1}{4}\right) w^2 d\xi \leq q(w)+  \|w\|_{L^2(\RR^+)}^2
$$
whence, using Cauchy-Schwarz inequality, 
\begin{align*}
\partial_{\tau} \|w\|_{L^2(\RR^+)}^2 + 2\left(1- \epsilon^{2\lambda} e^{-\lambda \tau} |\psi_\e| \right) q(w)
 \leq &  \frac{e^{\tau}}{\e^2}\partial_{yy} \|w\|_{L^2(\RR^+)}^2  \\
 &+ 2 \epsilon^{2\lambda} e^{-\lambda \tau} \left( (|\phi_\e| + |\psi_\e|+\frac12)\|w\|_{L^2(\RR^+)}^2 +\frac12 \|e^{\frac{\xi^2}{8}} f_\e\|_{L^2}^2 \right)
\end{align*}
If $\epsilon_0>0$ is small enough (depending on $\lambda$ and  $C_0$), $1- \epsilon^{2\lambda} e^{-\lambda \tau} |\psi_\e(\tau)|>0$ for any $\tau \geq 0$ and $\e  \in (0,\e_0)$, which combined with $q(w)\geq 0$ gives
\begin{equation}
\label{ineq sur norme de w}
\partial_{\tau} \|w\|_{L^2(\RR^+)}^2 \leq  \frac{e^{\tau}}{\e^2}\partial_{yy} \|w\|_{L^2(\RR^+)}^2 + \epsilon^{2\lambda} e^{-\lambda \tau}( C_1\|w\|_{L^2(\RR^+)}^2+C_2)
\end{equation}
where $C_1$ only depends on $\sup \{|\phi_\e(\tau)|, |\psi_\e(\tau)| \, , \, \tau \geq 0 \, , \, \e >0\}$ while $C_2$ depends on $f_\e$.
Set $h(\tau)$ solution to the ODE
$$
 \forall \tau \geq 0 \, , \quad h'(\tau)= \epsilon^{2\lambda}   e^{-\lambda\tau} (C_1 h(\tau)+C_2)  \, , \quad h(0)=\|w_0\|_{L^2(L^{\infty})}^2
$$
then $h$ is a supersolution to \eqref{ineq sur norme de w}  and for any $\tau \geq 0$, 
\begin{equation}
\label{edo}
h(\tau)=h(0)e^{\frac{C_1}{\lambda} \e^{2\lambda} (1-e^{-\lambda \tau })} + \frac{C_2}{C_1} \left( e^{\frac{C_1}{\lambda} \e^{2\lambda} (1-e^{-\lambda \tau })} -1 \right)
\end{equation}
If $\epsilon_0$ is small enough  (compared to $\lambda/C_1$), we can bound the second term as follows:
$$\|w(\tau)\|_{L^2(L^\infty)}^2 \leq C \left(\|w_0\|_{L^2(L^\infty)}^2+ \frac{C_2}{\lambda} \e^{2\lambda}\right)$$
This concludes the proof of lemma \ref{w bounded}. 
\rule{2mm}{2mm}

%%%%%%%%%%%%%%%%%%%%%%%%%%%%%%%
\noindent
\textbf{Proof of theorem \ref{lm 2.2} (continued).} We use the spectral property of ${\cal M}$ to decompose any solution $w$ to \eqref{dirichlet sym} as 
$$
\forall (\tau,\xi,y) \in \RR^+ \times\RR^+ \times \RR \, , \quad w(\tau,\xi,y)=\alpha(\tau,y) e_0(\xi) + r(\tau,\xi,y) 
$$
where $\alpha(\tau,y)=<w(\tau,\cdot,y),e_0>_{L^2(\RR^+)}$ so that $r$ is a transverse perturbation : for any $(\tau,y) \in \RR^+ \times \RR$, $<r(\tau, \cdot,y),e_0>_{L^2(\RR^+)}=0$. Projecting equation \eqref{dirichlet sym} on the null space of ${\cal M}$ gives
$$
\partial_{\tau} \alpha = \frac{e^{\tau}}{\e^2} \partial_{yy} \alpha + \e^{2\lambda} e^{-\lambda \tau} \left( (\phi_\e-\frac{\psi_\e}{\sqrt{\pi}}) \alpha -\psi_\e <r , e_0'+ \frac{\xi}{4} e_0>_{L^2(\RR^+)} + <e^{\frac{\xi^2}{8}} f_\e,e_0>_{L^2(\RR^+)} \right)
$$
while the equation satisfied by $r$ reads
\begin{equation}
\label{eq r}  
\partial_{\tau} r=  {\cal M}r + \frac{e^{\tau}}{\e^2} \partial_{yy} r  + \e^{2\lambda} e^{-\lambda \tau} \left( \phi_\e r +\psi_\e Q(\partial_{\xi} r -\frac{\xi}{4} r )+ \alpha \psi_\e \, Q (e_0'-\frac{\xi}{4} e_0) +Q(e^{\frac{\xi^2}{8}} f_\e) \right)
\end{equation}
where $P=1-Q$ is the projection onto the null space of ${\cal M}$.

\noindent
Since our idea is to find a dynamics similar to that of the heat equation, we introduce $\alpha_c$ solution to 
$$
\partial_{\tau} \alpha_c = \frac{e^{\tau}}{\e^2} \partial_{yy} \alpha_c \, , \quad \alpha_c(0,y)=\alpha(0,y)
$$
and denote $\beta=\alpha-\alpha_c$ the difference. Then, $\beta$ satifies $\beta(0,y)=0$ and 
\begin{equation}
\label{eq beta} 
\partial_{\tau} \beta = \frac{e^{\tau}}{\e^2} \partial_{yy} \beta +  \e^{2\lambda} e^{-\lambda \tau} \left((\phi_\e-\frac{\psi_\e}{\sqrt{\pi}})(\alpha_c+\beta) -\psi_\e <r, e_0'+ \frac{\xi}{4} e_0> +<e^{\frac{\xi^2}{8}} f_\e,e_0>\right)
\end{equation}
We shall prove that $\beta$ remains small for all time and that $r$ decreases exponentially fast to zero as time goes to infinity. 
Indeed, by the maximum principle and  lemma \ref{w bounded}, we get
\begin{align*}
\partial_{\tau} \beta \leq & \frac{e^{\tau}}{\e^2} \partial_{yy} \beta +  \e^{2\lambda} e^{-\lambda \tau}  \left( |\phi_\e| + \frac{\psi_\e}{\sqrt{\pi}} \right) |\beta| \\
& +  \e^{2\lambda} e^{-\lambda \tau}  \left( \left( |\phi_\e| + \frac{\psi_\e}{\sqrt{\pi}} \right) \|\alpha_c(0)\|_{L^{\infty}} + |\psi_\e|  \|e_0'+\frac{\xi}{4} e_0 \|_{L^2} \|r(\tau)\|_{L^2(L^{\infty})} + \|e^{\frac{\xi^2}{8}} f_\e\|_{L^2} \right)
\end{align*}
Define $h$ as a solution to the ODE
$$
h'(\tau)=\e^{2\lambda} e^{-\lambda \tau}(C_1 |h(\tau)| +C_2) \, , \quad h(0)=0
$$
where $C_1$ only depends on $\phi_\e$ and $\psi_\e$ while $C_2$ depends on $\phi_\e$, $\psi_e$, $f_\e$ and $\|w\|_{L^2(L^\infty)}$.
Then, $h$ is a supersolution to \eqref{eq beta} and dealing as in \eqref{edo}, we get for $\e_0$ small enough (compared to $\lambda/C_1$), 
\begin{equation}
\label{ineq de beta}
\forall \tau \geq 0 \, , \quad \|\beta(\tau)\|_{L^{\infty}(\RR)} \leq |h(\tau)| \leq e \frac{C_2}{\lambda} \e^{2\lambda}
\end{equation}
We shall now apply parabolic regularity to get the same bounds on the derivatives of $\beta$. For any $y_0 \in \RR$, set $\zeta=\e \, e^{-\frac{\tau}{2}} (y+y_0)$ and denote $B(\tau,\zeta)=B(\tau,\e e^{-\frac{\tau}{2}}(y+y_0))=\beta(\tau,y)$. Then, by \eqref{eq beta},
$$
\partial_\tau B = \partial_{\zeta\zeta}B + \frac{\zeta}{2} \partial_{\zeta}B + \epsilon^{2\lambda} e^{-\lambda \tau} \left((\phi_\e-\frac{\psi_\e}{\sqrt{\pi}})(\alpha_c+B)-\psi_\e <r,e_0'+\frac{\xi}{4} e_0>+<e^{\frac{\xi^2}{8}} f_\e,e_0>\right)
$$
The above bound on $\beta$ also gives $B$ uniformly bounded by $\epsilon^{2\lambda}$. Finally, the parabolic regularity applies for $|\zeta|<1$ and we get that the derivatives of $B$ are uniformly bounded by $\epsilon^{2\lambda}$. Coming back to $\beta$, we get the desired estimates since the bounds do not depend on $y_0$.

%%%%%%%%%%%%%%%%%%%%%%%%%
As far as $r$ is concerned, we compute an energy estimate to benefit from the spectral gap in self similar variables. Taking the $L^2$ scalar product of \eqref{eq r} with $r$ gives
\begin{align}
\label{ineq sur norme de r} & \partial_{\tau} \|r\|_{L^2(\RR^+)}^2 + 2 q(r)=  \frac{e^{\tau}}{\e^2} \left(\partial_{yy}\|r\|_{L^2(\RR^+)}^2 -  2 \|\partial_y r\|_{L^2(\RR^+)}^2\right) +2 \e^{2\lambda} e^{-\lambda \tau} \phi_\e \|r\|_{L^2(\RR^+)}^2\\
 \notag& \quad + 2\epsilon^{2\lambda} e^{-\lambda \tau}\left( \psi_\e<Q(\partial_\xi r- \frac{\xi}{4}r) + \alpha Q(e_0'-\frac{\xi}{4}e_0),r>_{L^2(\RR^+)} +<Q(e^{\frac{\xi^2}{8}} f_\e), r>\right) 
\end{align}
Since
$$
\left|<Q(\partial_{\xi}r-\frac{\xi}{4}r),r>_{L^2(\RR^+)}\right|  = \int_0^{\infty} \frac{\xi}{4} r^2 d\xi \leq \int_0^{\infty} \left(\frac{\xi^2}{16} + \frac{1}{4} \right) r^2 d\xi \leq q(r)+\|r\|^2_{L^2(\RR^+)}
$$
and
$$
\left| \alpha<Q(e_0'-\frac{\xi}{4}e_0),r>_{L^2(\RR^+)} \right| \leq \|\alpha(\tau)\|_{L^{\infty}} \|e_0'-\frac{\xi}{4}e_0\|_{L^2} \|r\|_{L^2} \leq \left(\|\alpha_c\|_{L^{\infty}} + \|\beta\|_{L^{\infty}}\right) \|r\|_{L^2}
$$
we get
\begin{align*}
\partial_{\tau} \|r \|_{L^2(\RR^+)}^2 & +  2(1-\epsilon^{2\lambda} e^{-\lambda \tau} |\psi_\e|) q(r) \leq  \frac{e^{\tau}}{\e^2} \partial_{yy}\|r\|_{L^2(\RR^+)}^2  + 2\epsilon^{2\lambda} e^{-\lambda \tau} ( |\phi_\e| +|\psi_\e|)\|r\|_{L^2}^2 \\
& + 2\epsilon^{2\lambda} e^{-\lambda \tau} \left( |\psi_\e| ( \|\alpha_c\|_{L^{\infty}} + \|\beta\|_{L^{\infty}}) \|r\|_{L^2} + \|e^{\frac{\xi^2}{8}} f_\e\|_{L^2} \|r\|_{L^2} \right)
\end{align*}
If $\epsilon_0$ is small enough (depending on $\lambda$ and  $C_0$), $1-  \epsilon^{2\lambda}  e^{-\lambda \tau} |\psi_\e| \geq \frac{3}{4}$ for any $\tau\geq 0$ and $\epsilon \in (0,\epsilon_0)$ which combined with $q(r)\geq \|r\|_{L^2}^2$, \eqref{ineq de beta} and lemma \ref{w bounded} gives
$$
\partial_{\tau} \|r\|_{L^2(\RR^+)}^2 + \frac{3}{2} \|r\|_{L^2}^2   \leq   \frac{e^{\tau}}{\e^2} \partial_{yy}\|r\|_{L^2(\RR^+)}^2  +   C\e^{2\lambda} e^{-\lambda \tau} 
$$
Define $h$ as a solution to the ODE
$$
h'(\tau)+\frac{3}{2} h(\tau)=  C \e^{2\lambda} e^{-\lambda \tau} \, , \quad h(0)=\|r_0\|_{L^2(L^{\infty})}^2
$$
Then, $h$ is a supersolution to \eqref{ineq sur norme de r} and
\begin{equation}
\label{ineq de r}
\forall \tau \geq 0 \, , \quad  \|r(\tau)\|_{L^2(L^{\infty})}^2 \leq h(\tau) \leq C \epsilon^{2\lambda} e^{-\lambda \tau} + e^{-\frac{3}{2}\tau}  \|r_0 \|_{L^2(L^\infty)}^2
 \end{equation}
We shall now apply parabolic regularity to get some bounds on $r$.
 For any $y_0 \in \RR$, set $\zeta=\epsilon \, e^{-\frac{\tau}{2}} (y+y_0)$ and denote $R(\tau, \xi,\zeta)=R(\tau,\xi,\e e^{-\frac{\tau}{2}}(y+y_0))=r(\tau,\xi,y)$. Then, by \eqref{eq r},
$$
\partial_\tau R = {\cal M} R +  \partial_{\zeta \zeta } R +\frac{\zeta}{2} \partial_\zeta R  +\epsilon^{2\lambda}  e^{-\lambda \tau} \, \left( \phi_\e R + \psi_\e Q(\partial_{\xi} R- \frac{\xi}{4} R) + \alpha \psi_\e Q(e_0'-\frac{\xi}{4} e_0) +Q(e^{\frac{\xi^2}{8}} f_\e) \right)
$$
 Moreover, by \eqref{ineq de r}, $\|R\|_{L^2(L^\infty)}^2 \leq C \epsilon^{2\lambda} e^{-\lambda \tau}$ and the parabolic regularity states that for any compact $K$ of $\RR^+$, there exists $C_K>0$ independent of $y_0$ such that for any $\tau >0$, $\xi \in K$ and $|\zeta|<1$, 
 $$
 \mbox{max } \left(|\partial_\tau R| \, , \, |\partial_\xi R| \, , \, |\partial_{\xi\xi} R| \, , \, |\partial_\zeta R| \, , \, |\partial_{\zeta\zeta} R| \right) \leq C_K  \e^\lambda  e^{-\frac{\lambda}{2} \tau}
 $$
 Coming back to $r$, we get 
 $$
 \mbox{max } \left(|\partial_\tau r| \, , \, |\partial_\xi r| \, , \, |\partial_{\xi \xi} r| \right) \leq C_K \e^\lambda e^{-\frac{\lambda}{2}\tau}
 $$
 while 
 $$|\partial_y r| \leq C_K  \e^{\lambda+1}  e^{-\frac{\lambda+1}{2}\tau} \, , \quad  |\partial_{yy} r| \leq C_K  \e^{\lambda+2} e^{-\frac{\lambda+2}{2}\tau}$$
 This implies the lemma with $\tilde{v}(\tau,\xi,y)=\frac{r(\tau,\xi,y)}{\xi} e^{-\frac{\xi^2}{8}} e^{\frac{\lambda}{2}\tau}$.
\rule{2mm}{2mm}

%%%%%%%%%%%%%%%%%%%%%%%%%%%%%%%%%%%%%%%%%%%%%%%%%%
\subsection{When the initial datum goes to 0 as $\vert y\vert$ goes to infinity}
%%%%%%%%%%%%%%%%%%%%%%%%%%%%%%%%%%%%%%%%%%%%%%%%%%%

\noindent
The result that we are going to prove is much simpler than Theorem \ref{lm 2.2}. We could use this last result, but we prefer to give a direct 
approach.
\begin{prop}
\label{p4.1} Let $v$ be a solution of \eqref{dirichlet}, with initial datum $v_0$ satisfying
\begin{enumerate}
\item $\sup_{y\in\RR}\Vert e^{\xi^2/8}v_0(.,y)\Vert_{L^2(\RR^+)}<+\infty$,
\item $\di\lim_{y\to\pm\infty}v_0(\xi,y)=0$, uniformly in $\xi \in\RR^+$.
\end{enumerate}
Then we have $v(\tau,\xi,y)=\xi\tilde v(\tau,\xi,y)$ with
$$
\lim_{\tau\to+\infty}\Vert \tilde v(\tau,.)\Vert_{L^\infty(\RR^+\times\RR)}=0.
$$
\end{prop}

\noindent
{\bf Proof.} Let us first make the following simplifying assumption: there is $A>0$ such that
\begin{equation}
\label{e4.3}
v_0(\xi,y)=0\quad\hbox{if $\vert y\vert\geq A$.}
\end{equation}
This allows us to pass to self-similar variables in $y$: 
$
\zeta=\e \di\frac{y}{\sqrt t}.
$
And so, \eqref{dirichlet} becomes
\begin{align}
\label{dirichlet1}
  \partial_{\tau} v = ({\cal L}+{\cal N})v + \e^{2\lambda} e^{-\lambda\tau} \left( \phi_\e(\tau)v+\psi_\e(\tau)\partial_{\xi} v +f_\e(\tau,\xi) \right) &\, ,  \quad \tau>0\, , \, \xi >0 \, , \,  \zeta \in \RR   \\
 \notag v(\tau,0,y)=0 &\, , \quad \tau >0 \, , \,  \xi =0 \, , \,  \zeta \in \RR
\end{align}
with ${\cal N}=\partial_{\zeta\zeta}+\di\frac12\zeta\partial_\zeta$. The spectrum of ${\cal N}$, in the space $L^2(\RR,e^{\zeta^2/8}d\zeta)$, is  $\{\di\frac{k}{2} \, , \, k \in \NN^*\}$. And so, writing $v(\tau,\xi,\zeta)=e^{-(\xi^2+\zeta^2)/8}w(\tau,\xi,\zeta)$ we obtain the following equation
for $w$:
\begin{align}
\label{dirichlet sym1}
 &\partial_{\tau} w = ({\cal M}+{\cal P}) w + \e^{2\lambda} e^{-\lambda \tau} \left((\phi_\e(\tau)-\frac{\xi}{4}\psi_\e(\tau))w+ \psi_\e(\tau)\partial_\xi w + e^{\frac{\xi^2+\zeta^2}{8}} f_\e(\tau,\xi)\right)  \\ 
\notag 
&w(\tau,0,y)=0  \quad  \tau>0  \, , \,  y \in \RR
\end{align}
where  ${\cal P} w= \partial_{\zeta\zeta} w +\left(\di\frac{1}{4}-\di\frac{\xi^2}{16}\right)  w$. We have, for all $w(\tau,\xi, \cdot)\in L^2(\RR)$:
$$
\int_\RR({\cal P}w)w \, d\zeta\geq\frac12\Vert w\Vert_{L^2(\RR)}^2.
$$
Arguing as in the proof of Theorem \ref{lm 2.2}, we obtain
\begin{equation}
\label{e4.4}
\Vert w(\tau,.)\Vert_{L^2(\RR^+\times\RR)}\leq e^{-\tau/2}\Vert w(0,.)\Vert_{L^2(\RR^+\times\RR)}.
\end{equation}
This proves the convergence to 0 of $v$. In order to suppress assumption \eqref{e4.3}, let us notice that, for all $\delta>0$, the function  $(v_0(\xi,y)-\delta)^+$ satisfies \eqref{e4.3}. Moreover, due to the convexity of $v\mapsto (v-\delta)^+$, the function $(v(\tau,\xi,\zeta)-\delta)^+$
is a sub-solution of \eqref{dirichlet sym1}. And so, we have $v(\tau,\xi,\zeta)\leq \overline v^\delta(\tau,\xi,\zeta)$ where $\overline v(\tau,\xi,\zeta)$
solves \eqref{dirichlet sym1} with initial datum $(v_0(\xi,y)-\delta)^+$. So $\overline v^\delta$  satisfies \eqref{e4.4}, which entails, by elliptic regularity, its convergence to 0 on every compact subset of $\RR_+\times\RR$. Because the zero-order coefficients of the equation \eqref{dirichlet1} are positive at infinity, the convergence holds in fact in $L^\infty(\RR^+\times\RR)$. By elliptic regularity, this is also true for $\partial_\xi v$. The mean value theorem implies the result. \rule{2mm}{2mm}

%%%%%%%%%%%%%%%%%%%%%%%%%%%%%%%%%%%%%%%%%%%%%%%%%
%%%%%%%%%%%%%%%%%%%%%%%%%%%%%%%%%%%%%%%%%%%%%%%%%
\section{General large time asymptotics for the full KPP equation, proof of theorem \ref{thm1}}
\label{s5}
%%%%%%%%%%%%%%%%%%%%%%%%%%%%%%%%%%%%%%%%%%%%%%%%%
%%%%%%%%%%%%%%%%%%%%%%%%%%%%%%%%%%%%%%%%%%%%%%%%%

\noindent 
Let $u_0 \in {\cal C}(\RR^2)$ satisfy assumption \eqref{assumption}, i.e. trapped between two translates of $1-H$.
 Denote $u$ the unique classical solution to \eqref{KPP} emanating from $u_0$ at time $t=1$. 
 
As announced in the introduction, we shall construct two functions $\bar{u}(t,x,y)$ and $\underline{u}(t,x,y)$, defined for 
$t>1$,  $\{x\leq t^\delta\}$ (with $\delta $ small to be chosen later) and $y\in\RR$, which will consist in solving
 equation \eqref{KPP} inside this region, with Dirichlet condition the trace, at  $\{x=t^\delta\}$, of a function which
  solves \eqref{KPP} approximately in the diffusive zone. We will see, in the next sections, that the functions
   $\bar{u}(t,x,y)$ and $\underline{u}(t,x,y)$ actually mimic the behaviour of the true solution $u(t,x,y)$. 

It will, however, be convenient to work in the self-similar coordinates. Let $w(\tau,\xi,\eta)$ be defined as in Section 
\ref{section equations}. Recall that $w$ satisfies \eqref{shift} with initial condition $w(0,\xi,\eta)=e^{\xi}u_0(\xi+2,y)$. 

  We will need the following frame, borrowed from \cite{NRR-Brezis}. Under the assumption \eqref{assumption}, there are functions $\eta_\pm(\tau)$ and $q_\pm(\tau)$, and constants $0<\eta_0<\eta_1$, depending only on $x_1$ and $x_2$, satisfying
$$\eta_0\leq \eta_-(\tau)\leq\eta_+(\tau)\leq\eta_1,\quad q_\pm(\tau)=O(e^{-\frac{\tau}{4}}),$$
and such that for any $\tau >0$, $\xi > \xi_\delta$, 
\begin{equation}
\label{e5.2}
\eta_-(\tau) \xi  e^{-\frac{\xi^2}{4}}-q_-(\tau)\xi e^{-\frac{\xi^2}{7}} \leq w(\tau,\xi,y)
\leq \eta_+(\tau)  \xi  e^{-\frac{\xi^2}{4}} +q_+(\tau)\xi e^{-\frac{\xi^2}{7}}  e^{-e^{\delta\tau}}
\end{equation}

To see it, it suffices to apply the paragraphs "An upper barrier"  in \cite{NRR-Brezis} to the solution of the 1D KPP 
equation emanating from $1-H(x-x_1)$ and "A lower barrier"' to that emanating from $1-H(x-x_2)$ and apply the comparison principle.

 In the sequel, for every small $\e>0$, we will set
\begin{equation}
\label{e5.10}
T_\e=\e^{-2} \mbox{ and } \tau_\e=\ln T_\e \mbox{ such that } \e=e^{-\frac{\tau_\e}{2}}.
\end{equation}
In the next two sections, we will seek to apply Theorem \ref{lm 2.2} with the initial datum
\begin{equation}
\label{e5.11}
w(\tau_\e,\xi,y) = e^\xi e^{\frac{T_\e}{2}} u(T_\e, \xi+2,y).
\end{equation}
Due to \eqref{e5.2}, we will be able to control this initial condition.

%%%%%%%%%%%%%%%%%%%%%%%%%%%%%%%%%%%%%%%%%%%%%%%%%%%%
\subsection{Diffusive supersolution}
%%%%%%%%%%%%%%%%%%%%%%%%%%%%%%%%%%%%%%%%%%%%%%%%%%%

\noindent
For any $\delta \in (0,\frac{1}{2})$, define $\xi_\delta=e^{-(\frac{1}{2}-\delta)\tau}$ which corresponds to $x=t^\delta$ in  self similar coordinates. Let $\bar{w}$ the solution to 
\begin{align}
\label{supersolution}
\partial_\tau \bar{w}={\cal L} \bar{w} + e^\tau \partial_{yy} \bar{w} - \frac{3}{2} e^{-\frac{\tau}{2} } \partial_\xi \bar{w} & \quad \tau \geq \tau_\e  \, , \,  \xi >-\xi_\delta \, , \, y \in \RR \\
\notag \bar{w}(\tau, \xi_\delta,y)= e^{-e^{\delta \tau}} &  \quad \tau \geq \tau_\e  \, , \,  \xi =-\xi_\delta \, , \, y \in \RR\\
\notag \bar{w}(\tau_\e, \xi, y)=w(\tau_\e,\xi,y) & \quad \tau = \tau_\e \, , \,  \xi  >-\xi_\delta \, , \, y \in \RR
\end{align}
 Then, $\bar{w}$ is a supersolution to \eqref{shift} for $\xi > -\xi_\delta$. Indeed, by definition \eqref{ref mvt} 
$$
w(\tau,\xi,y)=e^{-\frac{\tau}{2}+\xi e^{\frac{\tau}{2}}} u_1(e^{\tau}, \xi e^{\frac{\tau}{2}},y)
$$
$u_1$ being strictly uniformly bounded by $0$ and $1$, it follows that 
$$
\forall \tau \geq 0 \, , \, \forall y \in \RR \, , \quad 0< w(\tau,-\xi_\delta,y)<e^{-e^{\delta \tau}}
$$
while $\partial_\tau w(\tau,-\xi_\delta,y)=\left(\partial_t u_1 e^{\frac{\tau}{2}} - \frac{1}{2} (u_1+\partial_x u_1) e^{-(\frac12 -\delta)\tau} -\frac{1}{2} u_1 e^{-\frac{\tau}{2}} \right) e^{- e^{\delta \tau}}$ gives for $\delta>0$ small enough
$$
\forall \tau \geq 0 \, , \, \forall y \in \RR \, , \quad \left| \partial_\tau w(\tau,-\xi_\delta,y) \right| \leq C e^{- \delta e^{\delta \tau}}
$$
To simplify the moving Dirichlet boundary $\xi=\xi_\delta=e^{-(\frac{1}{2}-\delta)\tau}$, we now introduce a change of variables:
$$
\bar{w}(\tau,\xi,y)=\bar{p}(\tau-\tau_\e, \xi+\xi_\delta, y) + e^{- e^{\delta \tau}} \chi(\xi+\xi_\delta)
$$
where $\tau_\e$ is defined in \eqref{e5.10} and $\chi$ is a smooth monotonic function such that $\chi(\eta)=1$ for $\eta \in [0,1)$ and $\chi(\eta)=0$ for $\eta >2$. The function $\bar{p}(\tau',\eta,y)$ then satisfies (removing the primes) for any $\tau >0$, $\eta >0$ and $y \in \RR$, 
\begin{align}
\label{eq de p bar} \partial_\tau \bar{p}= {\cal L} \bar{p} + \frac{e^\tau}{\e^2} \partial_{yy} \bar{p} + \e^{1-2\delta} e^{-(\frac12 -\delta)\tau}  &\left(- \left( \delta + \frac{3}{2} \e^{2\delta} e^{-\delta \tau} \right) \partial_\eta \bar{p} + \Xi_\e(\tau,\eta) \right)  \\
\notag \bar{p}(\tau, 0, y)= 0 & \quad  \tau>0  \, , \,  \eta=0  \, , \, y \in \RR \\
\notag \bar{p}(0,\eta,y)=w(\tau_\e, \eta-\e^{1-2\delta},y)-e^{-1/\e^{2\delta}} \chi(\eta)& \quad\tau=0\, , \,  \eta >0 \, , \, y \in \RR 
\end{align}
where $\Xi_\e$ is a smooth function supported in $\eta \in [0,2]$ and uniformly bounded:
$$\exists C_\delta >0 \, | \, \forall \e >0 \, , \, \forall \tau \geq 0 \, , \, \forall \eta \geq 0 \, , \quad  |\Xi_\e (\tau,\eta)| \leq C_\delta$$
Choose $\lambda=\frac{1}{2}-\delta>0$,  $\phi_\e=0$, $\psi_\e(\tau)=-(\delta+\frac{3}{2}\e^{2\delta} e^{-\delta \tau})$ uniformly bounded in $\tau$ and $\e$ and $f_\e=\Xi_\e$ compactly supported in $\eta$ and uniformly bounded in $\tau$ and $\e$. Then, applying theorem \ref{lm 2.2}, we have for $\tau >\tau_\e$, $\xi >-\xi_\delta$, $y \in \RR$, 
$$
\bar{w}(\tau,\xi,y)= (\xi+\xi_\delta) \left( \frac{e^{-\frac{(\xi+\xi_\delta)^2}{4}}}{\sqrt{2\sqrt{\pi}}} \left(\bar{\alpha}_c(\tau-\tau_\e,y) +\bar{\beta}(\tau-\tau_\e,y)\right)+ e^{- \frac{\lambda}{2} (\tau-\tau_\e)} \tilde{p}(\tau-\tau_\e,\xi+\xi_\delta,y) \right)
$$
where for any $\tau >0$ and $y \in \RR$
$$
\partial_{\tau} \bar{\alpha}_c = \frac{e^{\tau}}{\e^2} \partial_{yy} \bar{\alpha}_c \, , \quad \bar{\alpha}_c(0,y)= \frac{1}{\sqrt{2\sqrt{\pi}}}  \int_0^{+\infty}\eta \left( w(\tau_\e, \eta-\e^{1-2\delta},y)-e^{-1/\e^{2\delta}} \chi(\eta) \right) d\eta
$$
and for any $\tau >0$
\begin{align*}
\|\bar{\beta}(\tau)\|_{L^\infty(\RR)} \leq C \epsilon^{1-2\delta} \, ,& \, \|\partial_\tau \bar{\beta} (\tau)\|_{L^\infty(\RR)} \leq C \epsilon^{1-2\delta} \\
\|\partial_y \bar{\beta} (\tau)\|_{L^\infty} \leq C \e^{2-2\delta} e^{-\frac{\tau}{2}} \, , & \,  \|\partial_{yy} \bar{\beta}(\tau)\|_{L^\infty} \leq C \e^{3-2\delta}  e^{-\tau}
\end{align*}
and for any $\tau >0$, $\xi \in K$ compact set of $\RR^+$, $y\in\RR$
\begin{align*}
\mbox{max }(|\tilde{p}(\tau,\xi,y)|, |\partial_\tau \tilde{p}| \, ,& \, |\partial_\xi \tilde{p} | \,  , \, |\partial_{\xi\xi} \tilde{p} | ) \leq C_K  \e^{\frac12-\delta}\\
|\partial_y \tilde{p}| \leq C_K \e^{\frac32-\delta} e^{-\frac{\tau}{2}} \, , &\quad   |\partial_{yy} \tilde{p}| \leq C_K \e^{\frac52-\delta} e^{-\tau} 
\end{align*}

%%%%%%%%%%%%%%%%%%%%%%%%%%%%%%%%%%%%%
\subsection{Diffusive subsolution}
%%%%%%%%%%%%%%%%%%%%%%%%%%%%%%%%%%%%%

\noindent
Since $0<w(\tau,\xi,y) \leq \bar{w}(\tau,\xi,y) \leq C(\xi+\xi_\delta)$ for some large $C>0$ and $\tau \geq \tau_\e$, the non linear term in \eqref{shift} can be bounded as follows: for any $\xi > \xi_\delta >e^{-\frac{\tau}{2}}$
$$
 e^{\frac{3}{2}\tau - \xi e^{\frac{\tau}{2}}} w^2 \leq C (\xi+\xi_\delta) e^{\frac32 \tau - \xi e^{\frac{\tau}{2}}} w  \leq 2 C e^{\frac32 \tau} \xi_\delta  e^{- \xi_\delta e^{\frac{\tau}{2}}} \leq 2 C e^{(1+\delta)\tau} e^{-e^{\delta \tau}}w \leq C_0 e^{-(\frac12 -\delta)\tau} w
$$
so that a subsolution to \eqref{shift} is given by
\begin{align}
\label{subsolution} 
 \partial_{\tau} \underline{w} = {\cal L} \underline{w} +  e^{\tau} \partial_{yy} \underline{w} -\frac{3 }{2} e^{-\frac{\tau}{2}}\partial_{\xi} \underline{w} +C_0 e^{-(\frac12-\delta)\tau} \underline{w} &\, ,  \quad  \tau>\tau_\e  \, , \,  \xi > \xi_{\delta} \, , \, y \in \RR \\
\notag
\underline{w}(\tau,\xi_{\delta},y)=0 & \, , \quad \tau >\tau_\e  \, , \,  \xi = \xi_{\delta} \, , \,  y \in \RR \\
\notag 
\underline{w}(\tau_\e,\xi,y)= w(\tau_\e,\xi,y) &\, , \quad \tau =\tau_\e  \, , \,   \xi > \xi_{\delta} \, , \, y \in \RR  \, , \,  
\end{align}

\noindent
As in the previous section, we simplify the moving Dirichlet boundary by defining $\eta=\xi-\xi_{\delta}$, $\tau'=\tau-\tau_\e$  and set $\underline{w}(\tau,\xi,y)=\underline{p}(\tau',\eta,y)=\underline{p}(\tau-\tau_\e,\xi-\xi_{\delta},y)$. Then, $\underline{p}$ satisfies (after dropping the primes) for any  $\tau>0$ , $\eta > 0$ and $y \in \RR$,  
\begin{align}
\label{subsolution shift} 
 \partial_{\tau} \underline{p} = {\cal L} \underline{p} +  \frac{e^{\tau}}{\e^2} \partial_{yy} \underline{p} +\e^{1-2\delta} e^{-(\frac12-\delta)\tau} &\left( C_0  \underline{p}+(\delta-\frac{3}{2} \e^{2\delta } e^{-\delta \tau}) \partial_{\eta} \underline{p}  \right)\\
 \notag
 \underline{p}(\tau,0,y)=0 & \, , \quad \tau \geq 0  \, , \,   \eta=0 \, , \, y \in \RR \\
\notag 
\underline{p}(0,\eta,y)= w(\tau_\e, \eta+\e^{1-2\delta},y) &\, , \quad \tau=0 \, , \, \eta > 0 \, , \, y \in \RR   
\end{align}
Choose $\lambda=\frac{1}{2}-\delta>0$, $\phi_\e=C_0$, $\psi_\e=\delta-\frac{3}{2} \e^{2\delta } e^{-\delta \tau}$ uniformly bounded in $\tau$ and $\e$ and  $f_\e=0$.
Then, applying theorem \ref{lm 2.2}, we have for $\tau >\tau_\epsilon$, $\xi>\xi_\delta$ and $y \in \RR$, 
$$
\underline{w}(\tau,\xi,y)= (\xi-\xi_\delta) \left( \frac{e^{-(\xi-\xi_\delta)^2/4}}{\sqrt{2\sqrt{\pi}}} \left(\underline{\alpha}_c(\tau-\tau_\epsilon,y) +\underline{\beta}(\tau-\tau_\e,y)\right)+ e^{- \frac{\lambda}{2} (\tau-\tau_\e)} \tilde{q}(\tau-\tau_\e,\xi-\xi_\delta,y) \right)
$$
where for any $\tau >0$ and $y \in \RR$, 
$$
\partial_{\tau} \underline{\alpha}_c = \frac{e^{\tau}}{\e^2} \partial_{yy} \underline{\alpha}_c \, , \quad \underline{\alpha}_c(0,y)= \frac{1}{\sqrt{2\sqrt{\pi}}}  \int_0^{+\infty}\eta\,  w(\tau_\epsilon,\eta + \e^{1-2\delta},y) d\eta
$$
and for any $\tau >0$
\begin{align*}
\|\underline{\beta}(\tau)\|_{L^\infty(\RR)} \leq C \epsilon^{1-2\delta} \, ,& \, \|\partial_\tau \underline{\beta} (\tau)\|_{L^\infty(\RR)} \leq C \epsilon^{1-2\delta} \\
\|\partial_y \underline{\beta} (\tau)\|_{L^\infty} \leq C \e^{2-2\delta} e^{-\frac{\tau}{2}} \, , & \,  \|\partial_{yy} \underline{\beta}(\tau)\|_{L^\infty} \leq C \e^{3-2\delta}  e^{-\tau}
\end{align*}
and for any $\tau >0$, $\xi \in K$ compact set of $\RR^+$, $y\in\RR$
\begin{align*}
\mbox{max }(|\tilde{q}(\tau,\xi,y)|, |\partial_\tau \tilde{q}| \, ,& \, |\partial_\xi \tilde{q} | \,  , \, |\partial_{\xi\xi} \tilde{q} | ) \leq C_K  \e^{\frac12-\delta}\\
|\partial_y \tilde{q}| \leq C_K \e^{\frac32-\delta} e^{-\frac{\tau}{2}} \, , &\quad   |\partial_{yy} \tilde{q}| \leq C_K \e^{\frac52-\delta} e^{-\tau} 
\end{align*}

%%%%%%%%%%%%%%%%%%%%%%%%%%%%%%%%%%%%%%%%%%%%%%%%%
\subsection{The proof of Theorem \ref{thm1}}
%%%%%%%%%%%%%%%%%%%%%%%%%%%%%%%%%%%%%%%%%%%%%%%%%

\noindent
It is now, just a matter of applying the preceding sections in the right order. Note that we have for any $\tau >\tau_\e$, $\xi >\xi_\delta$ and $y \in \RR$, 
$$
0\leq\overline w(\tau,\xi,y)-\underline w(\tau,\xi,y)\leq C\e^{1-2\delta}.
$$
Define $\overline u_+$ and $\underline u_+$ the function corresponding to $\overline w(\tau,0,y)$ and $\underline w(\tau,2 \xi_\delta,y)$ in the moving frame (see \eqref{ref mvt} to \eqref{var auto sim}):
$$
\overline u_+ (t,y)=e^{-t^\delta} t^{1/2}\overline w(\ln t, 0,y),
\quad \underline u_+ (t,y)=e^{-\delta} t^{1/2}\underline w(\ln t, 2t^{-(\frac12-\delta)},y).
$$
Both $\overline u_+$ and $\overline u_-$  have the form \eqref{e3.30}, with estimate \eqref{e3.10} and assumptions \eqref{hyp 1 sur a} and \eqref{hyp 2 sur a}. 
Indeed, (dealing for instance with $\overline u_+$, and the same holds for $\underline u_+$)
$$
\overline u_+(t,y)=t^{\delta} e^{-t^\delta-1/4t^{1-2\delta}} \, \frac{\overline a(t,y)+\overline b(t,y)}{\sqrt{2\sqrt{\pi}}}
$$
where $\overline a(t,y)= \overline \alpha_c(\ln(t\e^2),y)$ satisfies $\partial_t \overline a= \partial_{yy} \overline a$ for any $t>1$ with $\overline a(1,y)=\overline \alpha_c(0,y)$ and $|b(t,y)| \leq C(\e^{1-2\delta} +1/t^{\frac14-\delta/2})$. $\overline a$ satisfies \eqref{hyp 1 sur a} and \eqref{hyp 2 sur a} thanks to \eqref{e5.2}. 
Proposition \ref{p3.1} and theorem \ref{t3.1} therefore imply
\begin{align*}
& U_{c_*}(x- \ln (\underline{a}(t,y)-C\e^{1-2\delta}))-\frac{C}{\sqrt{t}}\leq u(t,x,y)\leq  U_{c_*}(x-\ln (\overline a(t,y)+C\e^{1-2\delta}))+\frac{C}{\sqrt{t}} \\
& |\overline a(t,y) - \underline a(t,y) | \leq C \e^{1-2\delta}
\end{align*}
Now we choose
$$a_0^\e(y)=\underline{a}(1,y)=\underline{\alpha}_c(0,y)=\frac{1}{\sqrt{2\sqrt{\pi}}}  \int_0^{+\infty}\eta\,  w(\tau_\epsilon,\eta + \e^{1-2\delta},y) d\eta,$$ 
this finishes the proof. \rule{2mm}{2mm}

%%%%%%%%%%%%%%%%%%%%%%%%%%%%%%%%%%%%%%%%%%%%%%%%%%%%%
%%%%%%%%%%%%%%%%%%%%%%%%%%%%%%%%%%%%%%%%%%%%%%%%%%%
\section{Examples of convergence and nonconvergence}
\label{gauche}
%%%%%%%%%%%%%%%%%%%%%%%%%%%%%%%%%%%%%%%%%%%%%%%%%%%%%%%%%%%%%%%%%%%%
%%%%%%%%%%%%%%%%%%%%%%%%%%%%%%%%%%%%%%%%%%%%%%%%%%%%%%%%%%%%%%%%%%%%

\noindent
This section is devoted to the consequences of Theorem \ref{thm1}, i.e the proof of theorem \ref{thm2}. We will first give an example of nonconvergence by exploiting the fact that some solutions of the heat equation do not converge to anything. In the next three sub-sections, we will give various cases of convergence: the simplest one is that of an initial datum tending, as $\vert y\vert\to\infty$, to a unique translate of $1-H$. The next one is when the initial datum tends to a $y$-periodic translate of $1-H$. The last one is when the initial datum tends to two different limits as $y\to\pm\infty$: here, we will still have convergence, but only on compact sets in $y$.

%%%%%%%%%%%%%%%%%%%%%%%%%%%%%%%%%%%%%%%%%%%%%
\subsection{Suitably oscillating initial data}
%%%%%%%%%%%%%%%%%%%%%%%%%%%%%%%%%%%%%%%%%%%%%

\noindent
The starting point of our construction is the following solution to the standard heat equation - see \cite{CE}, \cite{VZ}, where similar phenomena are discussed:
$$
\partial_t a=\partial_{yy} a,\quad\hbox{or, with the change of variables $\tau=\ln (t)$:}\quad \partial_\tau \alpha= e^{\tau}\partial_{yy}\alpha,
$$
with initial datum $\alpha(0,y)=a(1,y)=\alpha_M(y)$, $M>1$ will be chosen later. Consider two sequences $(t_n)_{n\in{\mathbb{N}}}$ and $(x_n)_{n\in{\mathbb{Z}}}$ satisfying the following five requirements:
\begin{enumerate}
\item $x_n=x_{-n}$ for $n\in\mathbb{N}$ 
\item the sequences $(t_n)_{n\in{\mathbb{N}}}$ and   $(x_n)_{n\in{\mathbb{N}}}$ are increasing,
\item $\di\lim_{n\to+\infty}\di\frac{x_{n+1}}{x_n}=+\infty,$
\item $\di\lim_{n\to+\infty}\di\frac{x_{n}^2}{t_n}=0,$ $\di\lim_{n\to+\infty}\di\frac{x_{n+1}^2}{t_n}=+\infty,$
\item For $n \in \NN$, $\alpha_M\equiv1$ on $(x_{2n},x_{2n+1})$, $\alpha_M\equiv M$ on $(x_{2n+1},x_{2n+2})$ and $\alpha_M$ even.
\end{enumerate}
An example is $t_n=\sqrt n(n!)$, $x_n=\sqrt{n!}$. We have then
\begin{equation}
\label{e6.1}
\lim_{n\to+\infty} a(t_{2n},0)=1,\quad \lim_{n\to+\infty} a(t_{2n+1},0)=M>1.
\end{equation}
Indeed,  we have for $t>1$ and $y \in \RR$
$$
a(t,y)=\frac1{\sqrt{4\pi (t-1)}}\int_{\RR}e^{-(y-y')^2/4(t-1)}\alpha_M(y')dy'=\frac1{\sqrt\pi}\int_{\RR}e^{-z^2}\alpha_M(y+2z\sqrt{t-1})dz,
$$
and so, because $\alpha_M$ is even, this reduces to
$$
a(t,0)=\frac2{\sqrt\pi}\int_0^{+\infty}e^{-z^2}\alpha_M(2z\sqrt{t-1})dz.
$$
Now, use the fact that $\alpha_M(y)=\bar\alpha_M \in\{1,M\}$ on $(x_n,x_{n+1})$:
\begin{align*}
a(t_n,0)=\frac{2}{\sqrt\pi} \bar\alpha_M \int_{x_n/(2\sqrt{t_n-1})}^{x_{n+1}/(2\sqrt{t_n-1})}e^{-z^2}dz& +\frac{2}{\sqrt\pi}\int^{x_n/(2\sqrt{t_n-1})}_0 e^{-z^2}\alpha_M(2z\sqrt{t_n-1})dz \\
&+\frac{2}{\sqrt\pi}\int_{x_{n+1}/(2\sqrt{t_n-1})}^{+\infty}e^{-z^2}\alpha_M(2z\sqrt{t_n-1})dz.
\end{align*}
Because of requirement 4. and the dominated convergence theorem, the last two terms go to 0 as $n\to+\infty$. And so we have
$$
a(t_n,0)=\frac{2\bar\alpha_M}{\sqrt\pi}\int_{x_n/(2\sqrt{t_n-1})}^{x_{n+1}/(2\sqrt{t_n-1})}e^{-z^2}dz+o_{n\to+\infty}(1)=\bar\alpha_M+o_{n\to+\infty}(1).
$$
This proves \eqref{e6.1}. Consider now the diffusive super and sub solutions $\overline w(\tau,\xi,y)$ and $\underline w(\tau,\xi,y)$ constructed in Section \ref{s5}, and respectively defined by \eqref{supersolution} and \eqref{subsolution}, with the common initial datum at time $\tau =0$
$$
\overline w(0,\xi,y)=\underline w(0,\xi,y)=\lambda\alpha_M(y)(1-H(\xi)),
$$
where $H$ is the Heaviside function, and $\lambda>0$ will be adjusted as the discussion proceeds. We have
$$
\biggl(\underline w(\tau,\xi,y),\overline w(\tau,\xi,y)\biggl)=\lambda\alpha(\tau,y)\biggl(\underline W(\tau,\xi),\overline W(\tau,\xi)\biggl)
$$
where $\overline W$ and $\underline W$ solve, respectively, \eqref{supersolution} and \eqref{subsolution} with no term $\partial_{yy}$. From Theorem \ref{lm 2.2},  there is 
$0<\underline\Lambda^\infty\leq\bar\Lambda^\infty$ such that
$$
\biggl(\underline W(\tau,\xi),\bar W(\tau,\xi)\biggl)\to_{\tau\to+\infty}\biggl(\underline\Lambda^\infty,\bar\Lambda^\infty\biggl)\xi e^{-\xi^2/4}.
$$
Notice $\underline \Lambda^\infty >0$ since we choose $\alpha_M \geq 1>0$. We choose $M>0$ large enough so that
$$
M\underline\Lambda^\infty>\bar\Lambda^\infty.
$$
And, finally, we choose $\lambda>0$ such that 
$$
u(1,x,y)=e^{-x} \lambda \alpha_M(y) (1-H(x))\leq 1-H(x).
$$
So we have
$$
\begin{array}{rll}
\bar u(t,1,0)=& e^{-1} \sqrt{t} \, \bar w(\ln t,1/\sqrt{t},0)\sim_{t\to+\infty}  \bar\Lambda^\infty e^{-1} \lambda a(t,0)\\
\underline u(t,1,0)=& e^{-1} \sqrt{t} \, \underline w(\ln t,1/\sqrt{t},0)\sim_{t\to+\infty}\underline\Lambda^\infty e^{-1} \lambda a(t,0)\\
\end{array}
$$
Because $\underline u(t,x,y)\leq u(t,x,y)\leq\bar u(t,x,y)$ we have, in the end:
$$
\limsup_{n\to+\infty}u(t_{2n},1,0)\leq\lambda e^{-1}\bar\Lambda^\infty,\quad \liminf_{n\to+\infty}u(t_{2n+1},1,0)\geq \lambda e^{-1}M\underline \Lambda^\infty.
$$
Thus,
$$
\liminf_{n\to+\infty}u(t_{2n+1},1,0)>\limsup_{n\to+\infty}u(t_{2n},1,0),
$$
which is our counterexample and proves theorem \ref{thm2}(1). \rule{2mm}{2mm}

%%%%%%%%%%%%%%%%%%%%%%%%%%%%%%%%%%%%%%%%
\subsection{Initial data tending to a limit}
%%%%%%%%%%%%%%%%%%%%%%%%%%%%%%%%%%%%%%%%

\noindent
Let us consider $u_0$ such that 
$$
\lim_{y\to\pm\infty}u_0(x,y)=u^+_0(x),
$$
uniformly with respect to $x\in\RR$. Recall that, for compatibility with \eqref{assumption}, we should have
$$
1-H(x-x_2)\leq u^+_0(x)\leq 1-H(x-x_1).
$$
Let $u^+(t,x)$ be the one-dimensional solution of \eqref{ref mvt} emanating from $u_0^+$ and 
$\sigma_\infty$ (see \eqref{alpha infini}) such that
$$
u^+(t,x)\longrightarrow_{t\to+\infty}U_{c_*}(x+\sigma_\infty).
$$
Standard arguments from the theory of semilinear parabolic equations yield
$$
\lim_{y\to\pm\infty}u(t,x,y)=u^+(t,x),
$$
uniformly in $x$ and locally uniformly in $t$. Let $w(\tau,\xi,y)$ be defined by
\eqref{var auto sim}, and $w^+(\tau,\xi)$ be the corresponding 1D solution. We still have
$$
\lim_{y\to\pm\infty}w(\tau,\xi,y)=w^+(\tau,\xi),
$$
uniformly in $\xi$ and locally uniformly in $\tau$. Consider 
$$
\tilde w(\tau,\xi,y)=w(\tau,\xi,y)-w^+(\tau,\xi),
$$
and $\e>0$. For $\tau\geq\tau_\e=-2{\mathrm{ln}}\e$, the function $\tilde w$ falls in
the assumptions of Proposition \ref{p3.1}. So, 
$$
\lim_{\tau \to+\infty}\tilde w(\tau,\xi,y)=0,
$$
uniformly in $\xi$ and $y$. This translates to $\tilde u(t,x,y)=u(t,x,y)-u^+(t,x)$.

%%%%%%%%%%%%%%%%%%%%%%%%%%%%%
\subsection{Initial data that are asymptotically periodic in $y$}
%%%%%%%%%%%%%%%%%%%%%%%%%%%%%

\noindent
Consider first an initial datum $u_0(x,y)$ that is periodic in $y$. The function 
$\alpha_c(\tau,y)$ defined in Theorem \ref{lm 2.2}  tends as $\tau\to+\infty$ to the average of
its initial datum. The $\omega$-limit set of $u_0$ for the full system \eqref{shift} is therefore made up of functions of the form 
$\alpha\xi^+e^{-\xi^2/4}$.
Because of the stability of these functions under the asymptotic equation of \eqref{shift},
the set $\omega(u_0)$ is made up of only one of these functions, say $\alpha_\infty\xi^+e^{-\xi^2/4}$.

\noindent
Let now be $u_0(x,y)$ and $u^+_0(x,y)$ such that 
$$
\lim_{y\to\pm\infty}\biggl(u_0(x,y)-u_0^+(x,y)\biggl)=0,\quad\hbox{uniformly in $x$.}
$$
Let $u^+(t,x,y)$ be the solution emanating from $u^+_0(x,y)$ and, as before,
$$
\tilde u(t,x,y)= u(t,x,y)-u^+(t,x,y).
$$
Arguing as in the preceding section, we obtain the uniform convergence of $\tilde u$ to 0
as $t\to+\infty$ and prove theorem \ref{thm2}(3). \rule{2mm}{2mm}

%%%%%%%%%%%%%%%%%%%%%%%%%%
\subsection{Initial data tending to two different limits}
%%%%%%%%%%%%%%%%%%%%%%%%%%%%

\noindent
Let us consider $u_0$ such that 
$$
\lim_{y\to+\infty}u_0(x,y)=u_0^+(x),\quad\lim_{y\to-\infty}u_0(x,y)=u_0^-(x),
$$
uniformly with respect to $x\in\RR$. Recall that, for compatibility with assumption \eqref{assumption}, we should have
$$
1-H(x-x_2)\leq u_0^+(x),u_0^-(x)\leq 1-H(x-x_1).
$$
Let us come back to equation \eqref{shift}. We use the self-similar variable $\zeta=\di\frac{y}{\sqrt t}$, and
discover that the function $\alpha_c(\tau,\zeta)$ tends, as $\tau\to+\infty$, to $\alpha_c^\infty$, the unique solution 
of
$$
\begin{array}{rll}
-\di\frac{d^2\alpha_c^\infty}{d\zeta^2}-\di\frac\zeta2\frac{d\alpha_c^\infty}{d\zeta}=&0,\quad\zeta\in\RR\\
\alpha_c(\pm\infty)=&e^{-\sigma_\infty^\pm}.
\end{array}
$$
We have $\alpha_c(0)=\di\frac{e^{-\sigma_\infty^+}+e^{-\sigma_\infty^-}}2$, which is the pointwise limit of 
$\alpha_c(\tau,y)$. And from Theorem \ref{lm 2.2}, we have
$$
\lim_{\tau\to+\infty}w(\tau,\xi,\zeta)=\alpha_c^\infty(\zeta)\xi^+e^{-\zeta^2/4}.
$$
Undoing this and reverting to $u$ proves theorem \ref{thm2}(2). \rule{2mm}{2mm}

%%%%%%%%%%%%%%%%%%%%%%%%%%%%%%%%%%%%%%%%%%%%%%%%%%%%%%%%%%%%%%%%%%%%
%%%%%%%%%%%%%%%%%%%%%%%%%%%%%%%%%%%%%%%%%%%%%%%%%%%%%%%%%%%%%%%%%%%%
%%%%%%%%%%%%%%%%%%%%%%%%%%%%%%%%%%%%%%%%%%%%%%%%%%%%%%%%%%%%%%
%%%%%%%%%%%%%%%%%%%%%%%%%%%%%%%%%%%%%%%%%%%%%%%%%%%%%%%%%%%%%%

{\footnotesize
%%%%%%%%%%%%%%%%%%%%%%%%%%%%%%%%%%%%%%%%%
%%%%%%%%%%%%%%%%%%%%%%%%%%%%%%%%%%%%%%%%%%
 
}

\end{document}